\documentclass{amsart}
\usepackage{paper_2_}

\title{On the Gorenstein property of local models for $GL_n$}
\author{Siyu Peng}
\address{Massachusetts Institute of Technology}
\email{pengs14@mit.edu}
\date{\today}

\begin{document}

\begin{abstract}
Local models are used to model, étale-locally, the singularities of integral models of Shimura varieties. They are defined in terms of linear algebra and admit explicit descriptions by matrix equations. We determine exactly which
local models for $\GL_n$ are Gorenstein. On one hand, this gives the analogous information for integral models of Shimura varieties for unramified unitary groups with parahoric level structure, and on the other hand also results in certain natural matrix equations defining non-Gorenstein singularities.
\end{abstract}

{\let\newpage\relax\maketitle}

\section{Introduction}

\subsection{Background}
Local models are certain schemes that have the same local behavior as certain integral models of Shimura varieties. In \cite{goertz}, G\"ortz uses the combinatorics of the local models for $\GL_n$ to produce an open covering of the special fiber by affine schemes defined via matrix equations with cyclic relations (see Section \ref{sec: obtaining eqns}). In this paper, we use these matrix equations to study the Gorenstein property. We shall prove that almost all local models are not Gorenstein. As a corollary, we are able to use properties of local models to prove certain results about matrix equations as well. For example: 

\begin{proposition} \label{ex: motivating example}
    Let $k$ be a field and $r, m > 1$. The affine scheme 
    \[
    \spec k[X_1, \dots, X_m] / (X_1 \dots X_m ,X_2 \dots X_m X_1 , \dots , X_m X_1 \dots X_{m-1} )
    \]
    for $r \times r$ matrices $X_i$ is not Gorenstein. 
\end{proposition}

Without the theory of local models, it is not clear how one would approach this. It is not even clear what the dimension of this scheme is, or if it is Cohen-Macaulay. With the correspondence given in \cite{goertz}, we may identify this as an open subscheme of the local model $M(rm, r(m-1), \{ 0, r, \dots, (m-1)r \})$ (defined in Definition \ref{def: gen local model} below), from which we may deduce that the dimension of this scheme is $r^2 (m-1)$, and by \cite{hsy26}, that it is Cohen-Macaulay.

In this example, we used information about the local model to learn something about matrix equations. In this paper, we shall mostly go in the opposite direction. In particular, we shall use these charts to translate the Gorenstein problem for local models into one with explicit equations. Afterwards, we may resolve it with commutative algebra. If one takes certain facts about the local model as a black box, the methods in this paper are elementary. 

Local models were first formalized by Rapoport-Zink \cite{RZ} in their study of integral models of Shimura varieties of PEL-type with parahoric level structure. A group-theoretic approach to its construction is given by Zhu \cite{zhu14} and extended by Pappas-Zhu \cite{pappas-zhu}. A helpful survey of the theory is \cite{PRS-survey}. 

We now give a preliminary definition of the local model for the general linear group $\GL_n$. The techniques presented in this paper will apply for the generalized local model (Definition \ref{def: gen local model}), where the conditions on the map $\phi_i$ are relaxed. 
\begin{definition}
    Let $O$ be a discrete valuation ring with uniformizer $\pi$ and fraction field $K$. Let $0 < r, m < n$ be integers and let $I = \{ i_0 < \dots < i_{m-1} \} \subseteq \{0, 1, \dots, n-1\}$ a nonempty subset. For $0 \leq j \leq m-1$, let $\phi_j: O^n \to O^n$ be diagonal maps with $\pi$ in the $i_j+1$ through $i_{j+1}$ (modulo $n$, with $i_{m} = i_0$) entries and $1$ in the other entries.  The \textbf{local model} $M(n, r, I)$ is the functor attaching to an $O$-algebra $S$ the set of commutative diagrams 
    \[
    \begin{tikzcd}
        F_0 \arrow[d] \arrow[r] & F_1 \arrow[r] \arrow[d] & \dots \arrow[r] & F_{m-1} \arrow[d] \arrow[r] & F_0 \arrow[d] \\
        S^n \arrow[r, "\phi_0"] & S^n \arrow[r, "\phi_1"] & \dots  \arrow[r] & S^n \arrow[r, "\phi_{m-1}"] & S^n 
    \end{tikzcd}
    \]
    where each $F_i$ is Zariski locally a direct summand of rank $r$. 
\end{definition}

The study of the singularities of local models is a central problem in this subject: In the $GL_n$ case, they étale-locally (and after base change to a suitable unramified extension of the reflex field) describe the singularities of the Rapoport-Zink integral models of Shimura varieties attached to unitary groups that are unramified at the prime $p$ under consideration, and with parahoric level structure at $p$. 

In the case of $\GL_n$, Görtz \cite{goertz} proved that the local model as defined by Rapoport and Zink is flat. Using the group-theoretic construction, Pappas-Zhu \cite{pappas-zhu} proved that the local model is normal and its special fiber is reduced. They also show that the irreducible components of the geometric special fiber are Cohen-Macaulay. When the parahoric level is special, the special fiber is irreducible and hence Cohen-Macaulay. He \cite{he13} proved the local model is Cohen-Macaulay for unramified groups, which is the case in $\GL_n$. This result was extended by Haines-Richarz \cite{hr22} and Fakhruddin-Haines-Lourenco-Richarz \cite{fh25} to tamely and wildly ramified groups with the exceptions when the residue characteristic is 2. He-Schremmer-Yu \cite{hsy26} proved the Cohen-Macaulay in general with a purely combinatorial approach.

\subsection{Gorenstein Local Models}

After establishing the Cohen-Macaulayness of local models, one may wonder if it is Gorenstein. A Cohen-Macaulay scheme that admits a dualizing complex, which is true in our setting, has a dualizing sheaf, and this sheaf is invertible exactly when the scheme is Gorenstein. Our result below gives a complete classification of the Gorenstein cases within the aforementioned class of Shimura varieties associated with the local model. 

This is interesting whenever one needs to work with the dualizing complex of Shimura varieties. For example, Fakhruddin-Pilloni \cite{Fakhruddin-Pilloni23} studied the integrality properties of Hecke operators on the
coherent cohomology of Shimura varieties. They utilize the Gorenstein property to construct a trace map for a morphism of Shimura varieties, which allows them to define such operators integrally. 

Let us collect the previously known results regarding the Gorenstein property here. Clearly, when the local model is regular, it is Gorenstein. For example, in the Drinfeld case (see \cite{r90}), the local models have semistable reduction and are regular. For a classification of these cases, see \cite{hpr20}. The case when the level is a special maximal parahoric was studied by Zhu in \cite{zhu14}. For unramified groups and when the level is Iwahori, ($I = \{ 0, \dots, n-1 \}$), He in \cite{he13} relates the special fiber of the local model with the closure of the Borel subgroup in the wonderful compactification. Combining this with the work of Brion-Polo \cite{bp00}, the Gorenstein property of these local models follows. For the ramified unitary local model of signature $(n-1, 1)$, the Gorenstein property of certain cases was investigated by Luo \cite{luo25}, and a determination of the Gorenstein property for the Bruhat-Tits strata of such local models is proved by He-Luo-Shi in \cite{hls26}. 

We collect the results for local models over $\GL_n$ in the following proposition. 

\begin{proposition}\label{prop: gorenstein cases}
    For each of the following cases, the local model $M(n, r, I)$ is Gorenstein:
    \begin{enumerate}
    \item $r = 1$ or $n-1$
    \item $|I| = 1$
    \item $I =\{ 0, 1, \dots, n-1\}$ 
    \item $I = \{ \nu, \nu+1 \}$ for any $\nu \in \{ 0, 1, \dots, n-1\}$ (when $\nu = n-1$, we set $\nu + 1 = 0$)
\end{enumerate}
For (1), (2), and (4), they are even regular. 
\end{proposition}
\begin{proof}

    \begin{enumerate}
        \item By Remark \ref{rem: duals}, we may assume that $r = n-1$. Then, it suffices to look at the open set obtained from Theorem \ref{thm: X_i's}:  
    \[
    U((t_i)_i) \cong k[x_0, \dots, x_{m-1}] / (x_{m-1} \cdots x_0 ). 
    \]
    We see that this is regular and hence Gorenstein.
    \item We see that $M(n, r, I)$ is by definition the standard Grassmannian $\mathrm{Gr}(r, n)$, which is smooth and thus Gorenstein.
    \item This is the Iwahori case, which follows from combining the work of He \cite{he13} to relate the local model to certain large Schubert varieties with the work of Brion-Polo \cite{bp00}, to conclude that they are Gorenstein. 
    \item This follows from a computation done by G\"ortz in Section 4.4.5 of \cite{goertz} 
    \end{enumerate}
\end{proof}

\subsection{Main Results}

In this paper, we shall prove that the remaining cases are not Gorenstein. 

\begin{theorem}[Main Theorem]
    For $n \geq 4$, $1 < r < n-1$, $1 < |I| < n-1$ except for $I = \{ \nu, \nu + 1 \}$, $M(n, r, I)$ is not Gorenstein.  
\end{theorem}

The sketch of the proof is as follows: we show first that the local models $M(4, 2, \{ 0, 2 \})$ (Corollary \ref{cor: M(4, 2, {0,2})}) and when $|I| = n-1$ (Corollary \ref{cor: |I| = n-1 is not Gorenstein}) are not Gorenstein. The former can be done by direct computation, whereas the latter is accomplished via induction, using that for a Gorenstein local ring, every quotient by an ideal generated by a regular sequence is again Gorenstein. (Here it is crucial that by flatness the dimensions of the special fibers of all local models are known.) We then choose certain open charts of the local model that reduce via the results Lemmas \ref{lem: remove 0's}, \ref{lem: remove 1's}, and Remark \ref{remark: remove rows} to these cases.

This paper begins with preliminary notions in Section \ref{sec: prelims}: the definition of the local model and how to translate charts into matrix equations in Section \ref{sec: obtaining eqns}. We also record some commutative algebra results that we shall use later in Section \ref{sec: gorenstein}. The proofs of our results will be given in Section \ref{sec: main results}. We begin with the base cases in Section \ref{sec: base cases}, which we shall use in the proof of the main theorem. We then take a detour and prove Proposition \ref{ex: motivating example} with the commutative algebra tools we developed. We conclude by proving the main theorem in Section \ref{sec: main thm}. 

\subsection{AI Disclosure} 

Anthropic's Claude Sonnet was used to help with literature review and editing the introduction, along with general typesetting. The rest of this work was derived without the use of generative AI. 

\subsection{Acknowledgements}
The author is deeply grateful for Professor Ulrich Görtz for proposing and advising this project. This project was funded by the MISTI-Germany program and the author would like to thank the University of Duisburg-Essen for their hospitality. 

\section{Preliminaries} \label{sec: prelims}

In this section, we shall give an explicit definition of the local model for $\mathrm{GL}_n$. In particular, we shall follow the work in Section 4 of \cite{goertz} to construct affine open subsets of the local model described by matrix equations. 

\subsection{Definitions}

Let $O$ be a discrete valuation ring with uniformizer $\pi$, residue class field $k$, and fraction field $K$. Fix some $n \geq 1$. For $i \in \{ 0, 1, \dots, n-1\}$, let $\Lambda_i$ be free $O$-modules of rank $n$ generated by $e^i_1, \dots, e^i_n$. 

\begin{definition}\label{def: gen local model}
    Let $m, n, r \in \ZZ_{>0}$ with $r < m \leq n$ and module homomorphisms $\phi_i: \Lambda_i \to \Lambda_{i+1}$ for $i = 0, \dots, m-1$ given by diagonal matrices with either $1$ or $\pi$ on the diagonal such that the composition $\phi_{m-1} \circ \phi_{m-2} \circ \dots \circ \phi_0$ is multiplication by $\pi$. The \textbf{generalized local model} $M(m, n, r, (\phi_i)_i)$ is the functor attaching to an $O$-algebra $S$ the set of commutative diagrams
    \[
    \begin{tikzcd}
        F_0 \arrow[d] \arrow[r] & F_1 \arrow[d] \arrow[r] & \dots  \arrow[r] & F_{m-1} \arrow[d] \arrow[r] & F_0 \arrow[d] \\ 
    \Lambda_{S, 0} \arrow[r, "\phi_0"] \ & \Lambda_{S, 1} \arrow[r, "\phi_1"]  & \dots \arrow[r] & \Lambda_{S, m-1} \arrow[r, "\phi_{m-1}"]  & \Lambda_{S, 0} \\
    \end{tikzcd}
    \]
    where $\Lambda_{S, i} =  \Lambda_i \otimes_O S$ is a free $S$-module of rank $n$ and each $F_i$ are Zariski-locally on $S$ a direct summand of $\Lambda_{S, i}$ of rank $r$.
\end{definition}

This functor is representable by a projective $O$-scheme. In particular, it is representable as a closed subscheme of a product of Grassmannians over $O$, $\prod_{j=0}^{m-1} \mathrm{Gr}(r, n)$. 

\begin{remark}\label{rem: duals} 
    We may extend the commutative diagram in Definition \ref{def: gen local model} to include the quotients: 
    \[
    \begin{tikzcd}
    0 \arrow[d] & 0 \arrow[d] & & 0 \arrow[d] & 0 \arrow[d] \\
    F_0 \arrow[d] \arrow[r] & F_1 \arrow[d] \arrow[r] & \dots  \arrow[r] & F_{m-1} \arrow[d] \arrow[r] & F_0 \arrow[d] \\ 
    \Lambda_{S,0} \arrow[r] \arrow[d] & \Lambda_{S, 1} \arrow[r] \arrow[d] & \dots \arrow[r] & \Lambda_{S, m-1} \arrow[r] \arrow[d] & \Lambda_{S, 0} \arrow[d] \\
    \Lambda_{S, 0}/F_0 \arrow[d] \arrow[r] & \Lambda_{S, 1} / F_1  \arrow[r] \arrow[d] & \dots \arrow[r] & \Lambda_{S, m-1} / F_{m-1} \arrow[r] \arrow[d] & \Lambda_{S, 0} / F_0 \arrow[d] \\
    0 & 0 &  & 0 & 0 
    \end{tikzcd}
    \]
    Since $F_i$ is locally a direct summand, each exact sequence is locally split and $\Lambda_{S, i} / F_i$ is also locally a direct summand of rank $n-r$. Since $\Lambda_{S, i} / F_i$ is finitely-presented, the sequence splits globally. The morphisms $\Lambda_{S, i} \to \Lambda_{S, i} / F_i$ induce a map on duals $(\Lambda_{S, i} / F_i)^\vee \to \Lambda_i^\vee$, and the maps $\Lambda_{S, i}/F_i \to \Lambda_{S, i+1}/F_{i+1}$ induces a map of duals in the other direction $(\Lambda_{S, i+1}/ F_{i+1})^\vee \to (\Lambda_{S, i} / F_i)^\vee$. Identifying $\Lambda_{S, i}^\vee \cong \Lambda_{S, i}$ and noting that $\phi_i^\vee = \phi_i$ as these are given by diagonal matrices, we have the diagram 
    \[
    \begin{tikzcd}
        (\Lambda_{S, 0} / F_0)^\vee \arrow[r] \arrow[d] & (\Lambda_{S, m-1} / F_{m-1})^\vee \arrow[r] \arrow[d] & \dots  \arrow[r] & (\Lambda_{S, 1} / F_1)^\vee \arrow[r] \arrow[d] & (\Lambda_{S, 0} / F_0)^\vee  \arrow[d] \\
        \Lambda_{S, 0} \arrow[r, "\phi_{m-1}"] & \Lambda_{S, m-1} \arrow[r, "\phi_{m-2}"] & \dots \arrow[r] & \Lambda_{S, 1} \arrow[r, "\phi_0"] & \Lambda_{S, 0} 
    \end{tikzcd}
    \]
    
    This defines an isomorphism 
    \[
    M(m, n, r, (\phi_i)_i) \cong M(m, n, n-r, (\phi_i)_{m-1-i}).
    \]
\end{remark}

\begin{definition}
    The \textbf{standard local model} denoted $M(n, r, I)$ is the generalized local model $M(m, n, r, (\phi_i)_i)$ with $m = |I|$ and maps $\phi_{j} = \mathrm{diag}(1^{i_j}, \pi^{i_{j+1} - i_j}, 1^{n-i_{j+1}})$ where $i_j$ is the $j$th element of $I$ and when $j = |I|$, we take $i_{j+1} := i_0 + n$ instead. 
\end{definition}

\begin{definition}
    The \textbf{special fiber} of the local model is the fiber over the residue field $k$. The \textbf{generic fiber} of the local model is the fiber over the fraction field $K$.
\end{definition}

\begin{remark} \label{rmk: special fiber}
    Over the generic fiber, $\pi$ is invertible and thus the maps $\phi_i$ become isomorphisms. Each point $(F_i)_i$ is then completely determined by the value of $F_0$, and so $(M(n, r, I))_K \cong \mathrm{Gr}(r, n)_K$, which is smooth. As the local model is flat \cite{goertz} and the generic fiber is smooth, it suffices in our study of the singularities of these models to consider only the special fiber. From this point onwards, we shall work in the special fiber, and by abusing notation, we shall write $\Lambda_i$ for $\Lambda_{k, i}$. 
\end{remark}

\subsection{Obtaining matrix equations} \label{sec: obtaining eqns}

Fix some local model $M(m, n, r, (\phi_i)_i)$. There exists an open covering by choosing values $t_i(j)$ satisfying the following condition:

\begin{condition}\label{condition} For $i \in \{ 0, \dots, m-1\}$ and $j \in \{1, \dots, n\}$ 
\begin{enumerate}
    \item $t_i(j) = 0$ or $1$ 
    \item For every $i$ and as $j$ ranges, $t_i(j)$ takes the value of $1$ exactly $r$ times 
    \item For fixed $j$ and as $i$ ranges, $t_i(j)$ is a cyclic permutation of $(1^l, 0^{m-l})$ for some $l$ which depends on $j$ 
    \item Whenever $t_i(j) = 1$ and $t_{i+1}(j) = 0$ where $i+1$ is taken modulo $m$, then $\phi_i(e^i_j) = 0$ 
\end{enumerate}
\end{condition}

\begin{remark}
    This is equivalent to Condition 1 in \cite{goertz}, where in Condition $1.i$, it should be $(t_0(j), \dots, t_{m-1}(j))$ instead of $(t_1(j), \dots, t_n(j))$. 
\end{remark}

It is convenient to organize the $t_i$'s in a table, which we shall refer to as the $t$-table.

\begin{example} \label{running example}
    Consider the standard local model $M(6, 2, \{0, 2, 3\})$. The maps $\phi_i$ are as follows: 
    \[
    \phi_0= \mathrm{diag}(0, 0, 1, 1, 1, 1), \quad \phi_1 = \mathrm{diag}(1, 1, 0, 1,1, 1), \quad \phi_2 = \mathrm{diag}(1, 1, 1, 0, 0, 0) .
    \]
    Below are two valid (as in, satisfying Condition \ref{condition}) $t$-tables: 
    \[
    \begin{tabular}{c|c|c|c|c|c|c}
        $i$ & 1 & 2 & 3 & 4 & 5 & 6 \\
        \hline 
        $t_0$ & 1 & 1 & 0 & 0 & 0 & 0 \\
        $t_1$ & 0 & 0 & 1 & 1 & 0 & 0 \\
        $t_2$ & 0 & 0 & 0 & 1 & 1 & 0
    \end{tabular} \quad \begin{tabular}{c|c|c|c|c|c|c}
         $i$ & 1 & 2 & 3 & 4 & 5 & 6 \\
         \hline 
        $t_0$ & 1 & 1 & 0 & 0 & 0 & 0 \\
        $t_1$ & 0 & 0 & 1 & 0 & 0 & 1 \\
        $t_2$ & 0 & 0 & 0 & 1 & 0 & 1
    \end{tabular}.
    \]
\end{example}

\begin{definition}
    Let $(t_i)_i$ satisfy Condition \ref{condition}. We define $U((t_i)_i)$ to be the open subset for which the quotient $\Lambda_i / F_i$ is a free module of rank $n-r$, generated by the $e^i_j$ such that $t_i(j) = 0$. 
\end{definition}

Fixing some $(t_i)_i$ satisfying Condition \ref{condition} and a point $(F_i)_i \in U((t_i)_i)$, letting the $e^i_j$'s for $t_i(j) = 0$ form a basis for $\Lambda_i / F_i$. We may view the maps $\Lambda_i / F_i \to \Lambda_{i+1} / F_{i+1}$ as $(n-r) \times (n-r)$ matrices $X_i$ with values in $S$. We aim to exhibit $U((t_i)_i)$ as an affine scheme corresponding to a ring of the form
\[
k[X_0, \dots, X_{m-1}] / (X_{m-1} X_{m-2} \dots X_0, \dots, X_0 X_{m-1} \dots X_1)
\]
where we write $k[X_0, \dots, X_{m-1}]$ to mean the polynomial ring generated by the entries of $X_0, \dots, X_{m-1}$. 

First, we may make two simplifications: if there is a column consisting entirely of 0's or entirely of 1's, we may remove it. 

\begin{lemma}[{\cite[Lemma 4.9]{goertz}}] \label{lem: remove 0's}
    Let $T_0 = \{ j \mid t_i(j) = 0 \, \forall i \}$. Let $u_i$ be the vector obtained by deleting the entries $j$ such that $j \in T_0$ from $t_i$. We have the isomorphism 
    \[
    U((t_i)_i) \cong U((u_i)_i) \times \mathbb{A}^{r|T_0|}
    \]
    where $U((u_i)_i) \subseteq M(m, n-|T_0|, r, (\psi_i)_i)$, the local model with chain $\overline{\Lambda_i} = \Lambda_i / \langle e^i_j \mid j \in T_0 \rangle$ and maps $\psi_i: \overline{\Lambda}_i \to \overline{\Lambda}_{i+1}$ induced by $\phi_i$. 
\end{lemma}

\begin{lemma}[{\cite[Lemma 4.10]{goertz}}] \label{lem: remove 1's}
    Let $T_1 = \{ j \mid t_i(j) = 1 \, \forall i \}$. Let $u_i$ be the vector obtained by deleting the entries $j$ such that $j \in T_1$ from $t_i$. We have the isomorphism 
    \[
    U((t_i)_i) \cong U((u_i)_i) \times \mathbb{A}^{(n-r)|T_1|} 
    \]
    where $U((u_i)_i) \subseteq M(m, n-|T_1|, r-|T_1|, (\psi_i)_i)$, the local model with chain $\overline{\Lambda_i} = \Lambda_i / \langle e^i_j \mid j \in T_1 \rangle$ and maps $\psi_i: \overline{\Lambda}_i \to \overline{\Lambda}_{i+1}$ induced by $\phi_i$. 
\end{lemma}

We refer the reader to the source for the proofs of these lemmas. 

\begin{remark} \label{remark: remove rows}
    Let $T_0, T_1$ be as above. Suppose there is some $i'$ such that $\phi_{i'}$ is multiplication by $0$ only at basis vectors $e^{i'}_j$ with $j \in T_0 \cup T_1$. Then, the associated $\psi_{i'}$ after applying the preceding two lemmas will be the identity, which can essentially be removed. In other words, we have the isomorphism 
    \[
    U((t_i)_i) \cong U((u_i)_{i \neq i'}) \times \AA^{r|T_0| + (n-r)|T_1|} 
    \]
    where $U((u_i)_i)\subseteq M(m-1, n-|T_0| - |T_1|, r - |T_1|, (\psi_i)_{i \neq i'})$ and $(u_i)_i$ is the vector obtained by deleting the entries $j$ such that $j \in T_0 \cup T_1$ from $t_i$.  
\end{remark}

In short, given a $t$-table, the above two lemmas and the previous remark allows us to remove columns that are all $0$'s or all $1$'s, along with rows $i$ for which the associated map $\phi_i$ is the identity, as the affine space factor does not affect the singularities of the chart. Let us assume that $(t_i)_i$ does not have a $j$ such that $t_i(j) = 0$ or $1$ for all $i$. 

\begin{notation}
    Suppose that there is a pair $i, j$ such that $t_i(j) = 0$. Let 
    \[
    \iota (i,j ) = \# \{ j' \mid j' \leq j \text{ and } t_i(j') = 0  \} .
    \]
    In other words, $\iota(i, j)$ is the value such that $t_i(j)$ is the $\iota(i, j)$th zero in its row in the $t$-table. 
\end{notation}

\begin{definition} \label{defn: X_i's} We define $(n-r) \times (n-r)$ matrices $X_i$ for each $i \in \{ 0, \dots, m-1 \}$ to be the set of matrices satisfying that for all $i, j$ be such that $t_i(j) = 0$ and $t_{i+1}(j) = 0$ (colloquially, there is a $0$ in the $(i, j)$ position with a $0$ below it in the $t$-table), $X_i$ has a $0$ in the $\iota(i, j)$ column except for the $\iota(i+1, j)$ entry, where it is 1. 
\end{definition}

\begin{remark}
    Note this is well-defined as for any fixed $i$, $\iota(i, j)$ is different for each valid $j$ such that $t_i(j) =0$, and thus, the $\iota(i, j)$ column of $X_i$ is set exactly once as we range over all pairs $(i, j)$. 
\end{remark}

\begin{example}
    For the open subset corresponding to the left $t$-table in Example \ref{running example}, we may remove the last column by Lemma \ref{lem: remove 0's}. Then, the associated $\mathfrak{X}_i$ is given by 
    \[
    X_0 = \begin{bmatrix}
        x^0_{11} & x^0_{12} &  0 \\
        x^0_{21} & x^0_{22} & 0 \\
        x^0_{31} & x^0_{32} & 1 
    \end{bmatrix}, \quad X_1 = \begin{bmatrix}
        1 & 0 & x^1_{11} \\
        0 & 1 & x^1_{21} \\
        0 & 0 & x^1_{31} 
    \end{bmatrix}, \quad X_2 = \begin{bmatrix}
        x^2_{11} & x^2_{12} & 1 \\
        x^2_{21} & x^2_{22} & 0 \\
        x^2_{31} & x^2_{32} & 0
    \end{bmatrix}
    \]
    for $x^i_{j \ell} \in S$.  
\end{example}

\begin{theorem}[{\cite[Theorem 4.12]{goertz}}]\label{thm: X_i's} Let $M(m, n, r, (\phi_i)_i)$ be a generalized local model. For some $(t_i)_i$ satisfying Condition \ref{condition} with no columns consisting entirely of 0's or entirely of 1's, the associated open subset $U((t_i)_i)$ is isomorphic to the affine scheme
    \[
\mathcal M = \spec k[x_{j\ell}^i]/(X_{m-1} X_{m-2} \cdots X_0 , X_{m-2} \cdots X_0X_{m-1}, \dots, X_0 X_{m-1} \cdots X_{1})
\]
where the $X_i$'s are given by Definition \ref{defn: X_i's} and the $x^i_{j \ell}$ are the nonconstant entries of these matrices. 
\end{theorem}

Again, we refer the reader to \cite{goertz} for details of the proof of this theorem. To give a sketch: we may define a morphism $U((t_i)_i) \to \mathcal{M}$ by mapping $(F_i)_i \mapsto (X_i)_i$, the matrices describing the maps $\Lambda_i / F_i \to \Lambda_{i+1} / F_{i+1}$ with respect to the basis of residue classes of $e^i_j$ such that ie. $t_i(j) = 0$. To build the maps in the other direction, let $\alpha_i: \Lambda_i \to \Lambda_i / F_i$ denote the projection maps onto these basis elements. This uniquely determines $F_i$ and $M_i: F_i \to \Lambda_i$ as the kernel of $\alpha_i$. In particular, $F_i$ is generated by $e^i_j - A^i_j$ where $t_i(j) = 1$ and $A^i_j$ is a linear combination in the $e^i_j$ such that $t_i(j) = 0$. Fixing this as a basis of $F_i$, we see that the $M_i$ will then be $n \times r$ matrices with $I_r$ in the rows $j$ for which $t_i(j) = 1$. Since the $r \times r$ matrices $N_i$ satisfy $\phi_i M_i = M_{i+1} N_i$ for each $i$, we notice that $N_i$ must be the rows of $\phi_i M_i$ which are the identity in $M_{i+1}$, ie. the rows $j$ such that $t_{i+1}(j) = 1$. 

\begin{equation} \label{super cool diagram}
    \begin{tikzcd}
    F_0 \arrow[d, "M_0"] \arrow[r, "N_0"] & F_1 \arrow[d, "M_1"] \arrow[r, "N_1"] & \dots  \arrow[r] & F_{m-1} \arrow[d, "M_{m-1}"] \arrow[r, "N_{m-1}"] & F_0 \arrow[d, "M_0"] \\ 
    \Lambda_0 \arrow[r, "\phi_0"] \arrow[d, "\alpha_0"] & \Lambda_1 \arrow[r, "\phi_1"] \arrow[d, "\alpha_1"] & \dots \arrow[r] & \Lambda_{m-1} \arrow[r, "\phi_{m-1}"] \arrow[d, "\alpha_{m-1}"] & \Lambda_0 \arrow[d, "\alpha_0"] \\
    \Lambda_0/F_0 \arrow[r, "X_0"] & \Lambda_1 / F_1  \arrow[r, "X_1"] & \dots \arrow[r] & \Lambda_{m-1} / F_{m-1} \arrow[r, "X_{m-1}"] & \Lambda_0 / F_0 
\end{tikzcd}
\end{equation}

\begin{example}
    Continuing with Example \ref{running example}, the $M_i$'s are given by  
    \[
    M_0 = \begin{bmatrix}
        1 & 0 \\
        0 & 1 \\
        a^0_1 & a^0_2 \\
        a^0_3 & a^0_4 \\
        a^0_5 & a^0_6 
    \end{bmatrix}, \quad M_1 = \begin{bmatrix}
        a^1_1 & a^1_2 \\
        a^1_3 & a^1_4 \\
        1 & 0 \\
        0 & 1 \\
        a^1_5 & a^1_6
    \end{bmatrix}, \quad M_2 = \begin{bmatrix}
        a^2_1 & a^2_2 \\
        a^2_3 & a^2_4 \\
        a^2_5 & a^2_6 \\
        1 & 0 \\
        0 & 1 
    \end{bmatrix} 
    \]
    where the $a^i_j$ are variables. The $N_i$'s are given by 
    \[
    N_0 = \begin{bmatrix}
        a^0_1 & a^0_2 \\
        a^0_3 & a^0_4
    \end{bmatrix}, \quad N_1 = \begin{bmatrix}
        0 & 1 \\
        a^1_1 & a^1_2
    \end{bmatrix}, \quad N_2 = \begin{bmatrix}
        a^2_1 & a^2_2 \\
        a^2_3 & a^2_4
    \end{bmatrix}.
    \]
\end{example}

We are able to give an equivalent description of $U((t_i)_i)$ with these $N_i$'s. Analogous to the $\iota(i, j)$, we have: 

\begin{notation}
    Define 
    \[
    \eta(i, j) = \# \{ j' \mid j' \leq j \text{ and } t_i(j') = 1\}.
    \]
    In other words, the value such that $t_i(j)$ is the $\eta(i, j)$th one in its row in the $t$-table.
\end{notation}

Notice the rows of $N_i$ that are rows of the identity matrix are exactly the $\eta(i+1, j)$ for which $t_i(j) = 1$ and $t_{i+1}(j) = 1$. Furthermore, the value of such a row would be the $\eta(i, j)$th row of the identity matrix. 

\begin{corollary} \label{cor: N_i's equations} For some $(t_i)_i$ satisfying Condition \ref{condition} and taking the $N_i$'s as above, $U((t_i)_i)$ is isomorphic to
    \[
    \spec k[N_0, \dots, N_{m-1}]/ (N_0 \cdots N_{m-1} = \dots = N_{m-1} N_0 \dots N_{m-2} = 0)
\]
\end{corollary}
Note that the order below is different from that of Theorem \ref{thm: X_i's}. 
\begin{proof}
    Let $X_i$ denote the family obtained via Definition \ref{defn: X_i's} from $(t_i)_i$ and $N_i$ the resulting map between $F_i$'s. 
    
    For $0 \leq i \leq m-1$, let $t_i' = t_{m-i-1}$ with the $0$'s and $1$'s switched, and let $(u_i)_i$ be given as $u_i = t_{i+1}'$, taken modulo $m-1$ (so $u_0 = t_{m-1}'$). In other words, the $t$-table associated to $(u_i)_i$ is the $t$-table of $(t_i)_i$ with the order of the rows reversed, the last row removed and inserted as the first row, and all 0's and 1's switched. Now consider $U((u_i)_i) \subseteq M(m, n, n-r, (\phi_{m-1-i})_i)$. Let $Y_i$ denote the result of applying Definition \ref{defn: X_i's} on $(u_i)_i$. We see that $Y_i = N_i^T$, as the definition of $N_i$ agrees with the description in Definition \ref{defn: X_i's} after we reverse the order of $\phi_i$; replace $i$ with $i-1$, which cyclically shifts the order of the maps by 1; and swap 0's and 1's, which corresponds to taking the transpose.

    Now it remains to show that $U((u_i)_i) \cong U((t_i)_i)$. By our discussion in Remark \ref{rem: duals}, we have the isomorphism 
    \[
    U((t_i)_i) \cong U((t_i')_i).
    \]
    Furthermore, it is clear by Theorem \ref{thm: X_i's} that cyclically shifting the maps $\phi_i$ will cyclically shift the $X_i$'s, which will yield the same relations. Thus, $U((t_i)_i) \cong U((t_i')_i) \cong U((u_i)_i)$. Finally, note
    \begin{align*}
        U((t_i)_i) &\cong k[N_0^T, \dots, N_{m-1}^T] / (N_{m-1}^T \cdots N_0^T= \dots = N_0^TN_{m-1}^T \dots N_{1}^T = 0) \\
        &\cong k[N_0, \dots, N_{m-1}] / (N_0 \cdots N_{m-1} = \dots = N_{m-1} N_0 \dots N_{m-2} = 0)
    \end{align*}
\end{proof}

\subsection{The Gorenstein Property}\label{sec: gorenstein}

As we will be working with affine schemes, we can check the Gorenstein property very explicitly with commutative algebra. 

\begin{definition} We take the definition of Gorenstein as in \cite[Chapter 21]{eisenbud}. Let $R$ be a Noetherian ring.

\begin{enumerate}
    \item Suppose $R$ is a local ring of dimension 0 with maximal ideal $\mathfrak{m}$. The \textit{socle} of $R$ is the annihilator of $\mathfrak m$. We say $R$ is \textit{Gorenstein} if its socle has dimension 1 as a vector space over $R/ \mathfrak m$. 
    \item Suppose $R$ is a local ring of dimension $d$ with maximal ideal $\mathfrak{m}$. $R$ is \textit{Gorenstein} if it is Cohen-Macaulay and for any regular sequence $f_1, \dots, f_d \in \mathfrak {m}$, $R / (f_1, \dots, f_d)$ is Gorenstein. 
    \item Suppose $R$ is any Noetherian ring. It is \textit{Gorenstein} if for all primes $P \subseteq R$, $R_P$ is Gorenstein. 
    \item A scheme $X$ is \textit{Gorenstein} if it is locally Noetherian and all its local rings are Gorenstein. 
\end{enumerate}
\end{definition}

When the dimension of the ring is small, it is not too hard to guess a regular sequence (as we do in the proof of Proposition \ref{prop: free 2x2 matrices}) and reduce to the 0-dimensional case, then explicitly compute the socle. However, it is difficult in general to find regular sequences. We shall use the following fact to determine when a sequence is regular. 

\begin{lemma}[{\cite[Theorem 17.4.iii]{matsumura}}] \label{lem: height and reg seq}
    Let $(R, \mathfrak m)$ be a Cohen-Macaulay local ring. A sequence $f_1, \dots, f_r \in \mathfrak m$ is regular if and only if the height of the ideal $(f_1, \dots, f_r)$ is $r$. 
\end{lemma}

With this, we prove the following lemma, which allows us to reduce determining Gorenstein singularities to smaller dimensions. 

\begin{lemma} \label{lem: reg seq}
    Let $X = \operatorname{Spec} R$ be an affine scheme that is of finite-type over a field $k$, equidimensional, and Cohen-Macaulay. Let $Z \subseteq X$ an equidimensional closed subscheme given by $Z = \operatorname{Spec}(R/(f_1, \dots, f_r))$ for $f_i \in R$ such that $\dim X - \dim Z = r$. Then, if $X$ is Gorenstein, $Z$ is Gorenstein as well. 
\end{lemma}

\begin{proof}
    Write $I = (f_1, \dots, f_r)$. Let $Q$ be a minimal prime over $I$. We shall show that the height of $Q$ is $r$. Since each minimal prime $Q$ corresponds to an irreducible component $C_Q \subseteq Z$, and $Z$ is equidimensional, we have $\dim C_Q= \dim Z$. By \cite[Theorem 12.2.9]{vakil}, 
    \[
    \operatorname{ht} Q = \dim X - \dim C_Q = \dim X - \dim Z  = r.
    \]
    This shows that the height of any minimal primes over $I$ is $r$. 

    Now suppose that $Z$ is not Gorenstein, so there is a point $p \in Z$ for which $\Oo_{Z, p}$ is not a Gorenstein local ring. Let $P \subseteq R$ be the prime associated to $p$, so $\Oo_{Z, p} = R_P$. Notice that 
    \[
    \mathrm{ht} (I_P) = \min_{Q \supseteq I_P} \dim R_Q
    \]
    In particular, any minimal prime over $I_P$ is a minimal prime over $I$ contained in $P$, which is a minimal prime over $I$. Thus, its height is $r$ and $(f_1, \dots, f_r)$ is regular in $R_P$ by Lemma \ref{lem: height and reg seq}. If $\Oo_{Z, p}$ is not Gorenstein, then $\Oo_{X, p}$ is not either. 
\end{proof}

To utilize this, we'll need to know the dimensions of these dense open sets of the local model. 

\begin{proposition} \label{prop: dim of local model}
    The special fiber of the local model $M(n, r, I)$ is Cohen-Macaulay and equidimensional of dimension $r(n-r)$ 
\end{proposition}

\begin{proof}
    The Cohen-Macaulay property for generalized local models is proved in \cite{he13}. 
    
    In \cite{goertz}, it was proven that the local model is flat over $O$, and thus the special fiber is equidimensional with the same dimension as that of the generic fiber. In the generic fiber, $\pi$ is invertible and thus the maps $\phi_i$ are isomorphisms. Each point $(F_i)_i$ is then determined by the value of $F_0$, and so $(M(n, r, I))_K \cong \mathrm{Gr}(r, n)_K$, which has dimension $r(n-r)$.
\end{proof}

\section{Main Results}\label{sec: main results}

In this section, we shall prove our main theorem. In \ref{sec: base cases}, we prove a specific class of local models are not Gorenstein. We shall do so by first explicitly exhibiting a base case that is not Gorenstein, then inductively applying Lemma \ref{lem: reg seq}. This technique may be used to prove Proposition \ref{ex: motivating example}, which we do in \ref{sec: proof of example}. Lastly, we prove the main theorem in \ref{sec: main thm}. This will be done by picking an appropriate open set $U((t_i)_i)$ that will reduce to one of the affine schemes given in \ref{sec: base cases}. 

\subsection{Base Cases} \label{sec: base cases}

In this section, we shall prove that the following base cases are not Gorenstein. 

\begin{proposition} \label{prop: free 2x2 matrices}
    The affine scheme $\spec k[X, Y]/(XY, YX)$ where $X, Y$ are free $2 \times 2$ matrices is not Gorenstein.
\end{proposition}

\begin{proof}
    We realize this as the open subset $U((t_i)_i)$ of the local model $M(4, 2, \{ 0, 2\})$ with $t$-table
    \[
    \begin{tabular}{c|c|c|c|c}
        $i$ & 1 & 2 & 3 & 4 \\
        \hline 
        $t_0$ & 1 & 1 & 0 & 0 \\
        $t_1$ & 0 & 0 & 1 & 1
    \end{tabular}.
    \]
    Write 
    \[
    X = \begin{pmatrix}
        x_1 & x_2 \\
        x_3 & x_4 
    \end{pmatrix}, \quad \quad Y = \begin{pmatrix}
        y_1 & y_2 \\
        y_3 & y_4
    \end{pmatrix}.
    \]
    Let us localize at the maximal ideal corresponding to the origin, $(x_1, x_2, x_3, x_4, y_1, y_2, y_3, y_4)$, and denote the resulting local ring as $R$. Consider the sequence
    \[
    (y_1 - x_1, y_2 - x_3, y_3 -x_2, y_4 -x_4-x_3) .
    \]
    Let $I \subseteq R$ be the ideal generated by this sequence. In $R/ I$, the relations $XY = YX = 0$ become 
    \begin{align*}
        0 = \begin{pmatrix}
        x_1 & x_2 \\
        x_3 & x_4
    \end{pmatrix} \begin{pmatrix}
        x_1 & x_3 \\
        x_2 & x_4 + x_3
    \end{pmatrix} &= \begin{pmatrix}
        x_1^2 + x_2^2 & x_1x_3 + x_2x_4 + x_2 x_3 \\
        x_1 x_3 + x_2 x_4 & x_3 ^2 + x_4^2 + x_3 x_4 
    \end{pmatrix}\\
    0 = \begin{pmatrix}
        x_1 & x_3 \\
        x_2 & x_4 + x_3
    \end{pmatrix}\begin{pmatrix}
        x_1 & x_2 \\
        x_3 & x_4
    \end{pmatrix} &= \begin{pmatrix}
        x_1^2 + x_3^2 & x_1 x_2 + x_3 x_4 \\
        x_1 x_2 + x_3 x_4 + x_3^2 & x_2^2 + x_4^2 + x_3 x_4
    \end{pmatrix} 
    \end{align*}
    Expanding and simplifying, we have
    \[
    x_1^2 = x_2^2 = x_3^2 = x_2 x_3 = 0, \quad x_4^2 = x_1 x_2 = - x_3 x_4, \quad x_1 x_3 = - x_2 x_4.
    \]
    Thus, $R/ I$ is given explicitly by $k[x_1, \dots, x_4, y_1, \dots, y_4]_{\mathfrak m_0}$ modulo these relations, where $\mathfrak m_0$ denotes the maximal ideal corresponding to the origin. One may check that $R/I$ is generated by the residue classes of monomials of degree at most 2 as a vector space over $k$. Thus, it has dimension 0. By Proposition \ref{prop: dim of local model}, the dimension of $R$ is $4$, by Lemma \ref{lem: height and reg seq}, the sequence generating $I$ is regular. We can manually check that the socle is generated by 3 elements (see below), and thus $R/I$ is not Gorenstein. 
    \[
    \begin{tikzpicture}[scale=0.9]
        \draw (1.5, 2) rectangle (2.5, 3); 
        \node at (2, 2.5) {$1$};
        \draw (0,1) rectangle (1,2);
        \draw (1,1) rectangle (2,2);
        \draw (2,1) rectangle (3,2);
        \draw (3,1) rectangle (4,2);
    
        \node at (0.5,1.5) {$x_1$};
        \node at (1.5,1.5) {$x_2$};
        \node at (2.5,1.5) {$x_3$};
        \node at (3.5,1.5) {$x_4$};
    
        \fill[gray!25] (0.5,0) rectangle (1.5,1);
        \fill[gray!25] (1.5,0) rectangle (2.5,1);
        \fill[gray!25] (2.5,0) rectangle (3.5,1); 
    
        \draw (0.5,0) rectangle (1.5,1);
        \draw (1.5,0) rectangle (2.5,1);
        \draw (2.5,0) rectangle (3.5,1); 
    
        \node at (1,0.5) {$x_1x_4$};
        \node at (2,0.5) {$x_2x_4$};
        \node at (3,0.5) {$x_4^2$};
    \end{tikzpicture}
    \]
\end{proof}

\begin{corollary} \label{cor: M(4, 2, {0,2})} 
    The local model $M(4, 2, \{ 0, 2\})$ is not Gorenstein. 
\end{corollary}

This concludes the first base case. The second is a family of base cases of the following form: 

\begin{proposition} \label{prop: base case 2}
    Consider the affine scheme
    \[
    U_m = \spec(k[n^i_{j\ell}]/(N_0 \cdots N_m, N_1 \cdots N_m N_0, \dots, N_m N_0 \cdots N_{m-1})
    \]
    where 
        \[
    N_i = \begin{pmatrix}
        0 & 1 \\
        n^i_{11} & n^i_{12} 
    \end{pmatrix} \quad \text{for } 0 \leq i \leq m-1 \text{ and} \quad N_m = \begin{pmatrix}
        n^m_{11}& n^m_{12} \\
        n^m_{21} & n^m_{22} 
    \end{pmatrix}
    \] 
    and $n^i_{j \ell}$ range over the entries of the matrices $N_i$ for $0 \leq i \leq m$. Then, $U_m$ is not Gorenstein for all $m \geq 2$. 
\end{proposition}

\begin{proof}
    First, we notice that this is the open subscheme of $M(m+2, 2, \{0, 1, \dots, m\})$ given by applying Corollary \ref{cor: N_i's equations} to the $t$-table 
    \[
        \begin{tabular}{c|c|c|c|c|c|c}
           $j$ & 1 & 2 & 3 & $\dots$& $m-1$ & $m$ \\
           \hline 
           $t_0$ & 1 & 1 & 0 & $\dots$ & 0 & 0 \\
           $t_1$ & 0 & 1 & 1 & $\dots$ & 0 & 0 \\
           $\vdots$ & $\vdots$ &  & $\ddots$ & $\vdots$ & $\vdots$ & \\
           $t_{m}$ & 0 & 0 & 0 & \dots & 1 & 1 
        \end{tabular}.
        \] 
    We shall proceed via induction on $m$. Consider the base case $m=2$. Let $R$ be the localization at the origin of $k[a_1, a_2, b_1, b_2, c_1\dots, c_4] / (N_0N_1N_2, N_1N_2 N_0, N_2 N_0 N_1)$ where 
    \[
    N_0 = \begin{pmatrix}
        0 & 1 \\
        a_1 & a_2
    \end{pmatrix}, \quad N_1 = \begin{pmatrix}
        0 & 1 \\
        b_1 & b_2 
    \end{pmatrix}, \quad N_2 = \begin{pmatrix}
        c_1 & c_2 \\
        c_3 & c_4
    \end{pmatrix}.
    \]
    Now consider the sequence 
    \[
    (a_1-c_2, a_2-c_3, b_1-c_1, b_2-c_4) 
    \]
    and let $I \subseteq R$ be the ideal generated by this sequence. With computer algebra, we compute that the relations of $R / I$ are given by 
    \[
    k[c_1, c_2, c_4]_{(c_1, c_2, c_4)} / (c_1^2, c_2c_4, c_2^2 + c_1 c_4,  c_4^2 + c_1c_2,  c_1c_4^2,  c_4^3)
    \]
    One checks that $R/I$ is a finite-dimensional vector space over $k$ with generators
    \[
    \begin{tikzpicture}[scale=0.9]
    \draw (1, 2) rectangle (2, 3); 
    \node at (1.5, 2.5) {$1$};
    \draw (0,1) rectangle (1,2);
    \draw (1,1) rectangle (2,2);
    \draw (2,1) rectangle (3,2);

    \node at (0.5,1.5) {$c_1$};
    \node at (1.5,1.5) {$c_2$};
    \node at (2.5,1.5) {$c_4$};

    \fill[gray!25] (0.5,0) rectangle (1.5,1);
    \fill[gray!25] (1.5,0) rectangle (2.5,1);

    \draw (0.5,0) rectangle (1.5,1);
    \draw (1.5,0) rectangle (2.5,1);

    \node at (1,0.5) {$c_1c_2$};
    \node at (2,0.5) {$c_1c_4$};
\end{tikzpicture}.
\]
    It is hence 0-dimensional. Furthermore, its socle is generated by the two elements $c_1c_2$ and $c_1c_4$. By Proposition \ref{prop: dim of local model}, $\dim R = 4$, so we may apply Lemma \ref{lem: height and reg seq} to conclude that $U_2$ is not Gorenstein. 
    
    Now for $m > 2$, we shall inductively embed $U_{m-1}$ into $U_m$. Note that $\dim U_m - \dim U_{m-1} = 2$. We claim that we may embed the former into the latter by setting 
    \[
    N_{m-1} = \begin{pmatrix}
        0 & 1 \\
        1 & 0
    \end{pmatrix},
    \]
    ie. by quotienting out by the sequence $(n^{m-1}_{11} - 1,n^{m-1}_{12} )$, working locally at the point where $n^{m-1}_{11} -1$ and the rest are 0. For brevity, let us denote this permutation matrix by $\omega$. We observe that if we replace $N_m$ with $\omega N_m$, we have the substitution
    \[
    N_0 \cdots \omega  N_{m}, \dots, \omega N_{m} \cdots N_{m-2} \rightsquigarrow N_0 \dots N_{m-2} N_{m}, \dots , N_{m} N_0 \cdots N_{m-2} 
    \]
    which are exactly the relations defining $U_{m-1}$. As for the last relation $N_{m} N_0 \dots N_{m-1}=0$, we notice that making this substitution yields 
    \[
    \omega N_{m} N_0 \dots N_{m-2} \omega  = 0 \iff N_{m} N_0 \dots N_{m-2}  = 0
    \]
    and so does not contribute any more relations. We may then apply Lemma \ref{lem: reg seq} to conclude that $U_m$ is not Gorenstein as $U_{m-1}$ is not. 
\end{proof}

\begin{corollary} \label{cor: |I| = n-1 is not Gorenstein}
    The local model $M(m, 2, \{0, 1, \dots, \hat i, \dots, m-1\})$, where $\hat i$ denotes that the $i$ is omitted, is not Gorenstein. 
\end{corollary}

\begin{proof}
    We may obtain the equations in Proposition \ref{prop: base case 2} by applying Corollary \ref{cor: N_i's equations} to the $t$-table 
    
    \[
        \begin{tabular}{c|c|c|c|c|c|c|c|c|c|c|c}
           $j$ & 1 & 2 & 3 & $\dots$& $i$ & $i+1$ & $i+2$ & $i+3$ & $\dots$ & $m+1$ & $m+2$ \\
           \hline 
           $t_0$ & 1 & 1 & & & & & & & & \\
           $t_1$ &  & 1 & 1 & & & & & & & \\
           $\vdots$ & & & $\ddots$ & $\ddots$ & & & & & &  \\
           $t_{i-1}$ & & & & & 1 & 1 & & & &  \\
           $t_i$ & & & & & & &  1 & 1  & &  \\
           $\vdots$ & & & & & & & & $\ddots$ & $\ddots$& & \\
           $t_{m-1}$ & & & & & &  & & &  & 1 & 1  \\
           $t_{m}$ & 1 & & & & & & & & & & 1
        \end{tabular}.
    \]
    In words, it should form a staircase pattern, where the second $1$ in a row lies above the first $1$ in the next row, with the exception of one large step at exactly one spot. 
\end{proof}

\subsection{Proof of Proposition \ref{ex: motivating example}} \label{sec: proof of example}

In this section, we shall prove Proposition \ref{ex: motivating example}. Though this comes as a corollary of the main theorem \ref{thm: main}, we may prove this directly by starting with the variety in Proposition \ref{prop: free 2x2 matrices}, increasing the dimensions of the matrices $X$ and $Y$ in Proposition \ref{prop: I = 0, n}, then increasing the number of such free matrices. The proofs of this section illustrate more methods of applying Lemma \ref{lem: reg seq}. 

\begin{proposition} \label{prop: I = 0, n}
    For $n \geq 2$, the local model $M(2n, n, I=\{0, n\})$ is not Gorenstein. 
\end{proposition}
\begin{proof}
    We shall proceed via induction on $n$. The base case of $n=2$ was completed in Proposition \ref{prop: free 2x2 matrices} above. For a fixed $n$, we shall consider the following open subset $U_{2n, n, \{0, n\}} :=U((t_i)_i) \subseteq M(2n, n, \{ 0, n\})$ where the $t_i$ are 
    \[
    \begin{tabular}{c|c|c|c|c|c|c}
        $i$ & 1 & \dots & $n$ & $n+1$ & \dots & $2n$ \\
        \hline 
        $t_0$ & 1 & \dots & 1 & 0 & \dots & 0 \\
        $t_1$ & 0 & \dots & 0 & 1 & \dots & 1 
    \end{tabular}.
    \]
    Our goal is to embed $U_{2n, n, \{0, n\}} \hookrightarrow U_{2(n+1), n+1, \{0, n+1 \} }$ as a closed subscheme cut out by a regular sequence. The former is given by $\mathrm{Spec}(k[X, Y]/(XY, YX))$ where $X, Y$ are $n \times n$ matrices and the latter is given by $\mathrm{Spec} (k[X', Y']/(X'Y', Y'X'))$ where $X', Y'$ are $(n+1) \times (n+1)$ matrices. Consider viewing $X, Y$ as the top left $n \times n$ submatrix of $X'$ and $Y'$ respectively: 
    \[
    X'= \begin{pNiceArray}{ccc|c}
        \Block{3-3} <\Large> X & & & x_{n+1} \\
        & & & \vdots \\
        & & & x_{2n} \\
        \hline
        x_1 & \dots & x_n & x_0 
    \end{pNiceArray}, \quad Y'= \begin{pNiceArray}{ccc|c}
        \Block{3-3} <\Large> Y & & & y_{n+1} \\
        & & & \vdots \\
        & & & y_{2n} \\
        \hline
        y_1 & \dots & y_n & y_0 
    \end{pNiceArray}.
    \]
    Setting $x_0 = 1$ and $x_1, \dots, x_{2n},y_0,  y_1, \dots, y_{2n} = 0$ gives us a closed embedding of $U_{n, 2n, \{0, n\}}$ into $U_{n+1, 2(n+1), \{ 0, n+1\}}$. However, $(x_0-1, x_1, \dots, x_{2n}, y_0, \dots, y_{2n})$ is not a regular sequence. Instead, we claim that $(x_0-1, x_1, \dots, x_{2n})$ is regular and $V(x_0-1, x_1, \dots, x_{2n}) = V(x_0-1, x_1, \dots, x_{2n}, y_0, \dots, y_n)$ as subvarieties of $U_{2(n+1), n+1, I=\{ 0, n+1 \}}$. 

    One inclusion is clear, and so for equality, it remains to show that if we impose $x_0=1, x_1, \dots, x_{2n} = 0$, we will have $y_i = 0$ for all $0 \leq i \leq 2n$. Consider
    \[
     \begin{pNiceArray}{ccc|c}
        \Block{3-3} <\Large> X & & & 0 \\
        & & & \vdots \\
        & & & 0 \\
        \hline
        0 & \dots & 0 & 1 
    \end{pNiceArray} \begin{pNiceArray}{ccc|c}
        \Block{3-3} <\Large> Y & & & y_{n+1} \\
        & & & \vdots \\
        & & & y_{2n} \\
        \hline
        y_1 & \dots & y_n & y_0 
    \end{pNiceArray} = 0 .
    \]
    Taking the last row of $X'$ and each column of $Y'$, we see that $y_0, \dots, y_n = 0$. Similarly, the equations of $Y'X' = 0$ yield the relations $y_{n+1} , \dots, y_{2n} = 0$.

    The difference in dimensions of $M(n, 2n, \{ 0 , n\})$ and $M(n+1, 2(n+1), \{ 0, n+1\})$ is $(n+1)^2-n^2 = 2n+1$, which is the length of the proposed sequence. We may then apply Lemma \ref{lem: reg seq} to conclude. 
\end{proof}

Now we may increase the number of free matrices to conclude. Note the reindexing of the matrices. 

\begin{proposition}[Proposition \ref{ex: motivating example}] 
    For $r > 1, m> 1$, the affine scheme 
    \[
    \spec k[X_0, \dots, X_{m-1}] / (X_{m-1} \cdots X_0 ,X_{m-2} \cdots X_{0} X_{m-1} , \dots , X_0 X_{m-1} \cdots X_{1})
    \]
    for $r \times r$ matrices $X_i$ is not Gorenstein. 
\end{proposition}

\begin{proof}
    Let 
    \[
    R =k[X_0, \dots, X_{m-1}] /  (X_{m-1} \cdots X_0 ,X_{m-2} \cdots X_{0} X_{m-1} , \dots , X_0 X_{m-1} \cdots X_{1}).
    \]
    Notice that $\spec R$ corresponds to the open subset $U((t_i)_i)$ of the local model $M(rm, r(m-1), \{0, r, 2r, \dots, (m-1)r \})$ where $(t_i)_i$ are given by 
    \begin{table}[h]
        \centering
        \small
        \begin{tabular}{c|c|c|c|c|c}
        $i$ & $1, \dots, r$ & $r+1, \dots, 2r$ & \dots & $(m-2)r + 1, \dots, (m-1)r$ & $(m-1)r+1, \dots, mr$ \\
        \hline 
        $t_0$ & $1^{(r)}$ & $1^{(r)}$ & \dots & $1^{(r)}$ & $0^{(r)}$ \\
        $t_1$ & $0^{(r)}$ & $1^{(r)}$ & \dots & $1^{(r)}$ & $1^{(r)}$ \\
        $t_2$ & $1^{(r)}$ & $0^{(r)}$ & $\ddots$ & $1^{(r)}$ & $1^{(r)}$ \\
        $\vdots$ & \vdots & $\ddots$ & $\ddots$ & $\vdots$ & $\vdots$  \\ 
        $t_{m-1}$ & $1^{(r)}$ & $1^{(r)}$ & $\ddots$ & $0^{(r)}$ & $1^{(r)}$
    \end{tabular}
    \end{table}

    Thus, the dimension of $R$ is $r^2(m-1)$. We would like to embed the open subset $U_{2r, r, \{ 0, r\}} = \spec k[X_0, X_1] / (X_1 X_0, X_0 X_1)$ from Proposition \ref{prop: I = 0, n} into $\spec R$. We notice that $U_{2r, r, \{ 0, r\}}$ is cut out by the equations $X_i - I$ for $i > 2$, where $I$ is the identity matrix. This yields $(m-2)r^2$ equations, which is exactly the difference in dimensions between the two spaces and thus we may apply Lemma \ref{lem: reg seq}. 
\end{proof}

\subsection{Proof of Main Theorem} \label{sec: main thm}

\begin{theorem} \label{thm: main}
    For integers $r < n$ and $I = \{ i_0, \dots, i_{m-1}\}  \subseteq \{ 0, 1, \dots, n-1\}$ satisfying all of the following conditions: \begin{itemize}
        \item $n \geq 4$, 
        \item $1 < r, m < n-1$, and 
        \item $I \neq \{ \nu, \nu + 1 \}$ for some $\nu \in \{ 0, \dots, n-1\}$, 
    \end{itemize}
    The local model $M(n, r, I)$ is not Gorenstein.
\end{theorem}

\begin{proof}
    By picking an appropriate choice of $(t_i)_i$ described below, one can apply Lemmas \ref{lem: remove 0's} and \ref{lem: remove 1's}, along with Remark \ref{remark: remove rows} to the associated open set $U((t_i)_i)$ to show that $U((t_i)_i)$ is isomorphic to the product of affine space with one of the following affine varieties: 
    \begin{itemize}
        \item (\ref{prop: free 2x2 matrices}) $\spec(k[X, Y] / (XY, YX))$ for $2 \times 2$ matrices with free entries $X, Y$. 
        \item (\ref{prop: base case 2}) $\spec(k[n^i_{j \ell}]/(N_0 \cdots N_m, N_1 \cdots N_m N_0, \dots, N_m N_0 \cdots N_{m-1})$ \\ where 
        \[
    N_i = \begin{pmatrix}
        0 & 1 \\
        n^i_{11} & n^i_{12} 
    \end{pmatrix} \quad \text{for } 0 \leq i \leq m-1 \text{ and} \quad N_m = \begin{pmatrix}
        n^m_{11}& n^m_{12} \\
        n^m_{21} & n^m_{22} 
    \end{pmatrix}
    \] 
    \end{itemize}
    which we proved in Section \ref{sec: base cases} are non-Gorenstein. Thus, $U((t_i)_i)$ and the associated local model will not be Gorenstein.  
    
    By Remark \ref{rem: duals}, we may assume that $r \leq \frac{n}{2}$. Define $(t_i)_i$ as follows:  
    \begin{enumerate}
        \item If $m + r-1 \geq n$, then note that $m \geq n - r + 1  \geq \frac{n}{2}+1$. Thus, there is some $\mu \in I$ such that $\phi_{\mu}(e^{\mu}_{\nu}) = 0$ for exactly one $\nu$, and for all other $j$, $\phi_{\mu}(e^{\mu}_j) = e^{\mu+1}_j$. Set $t_i(\nu) = 1$ for all $i$. To fill this row such that it satisfies Condition \ref{condition}, we may set it to be the same as the row below it after following the rest of the procedure. By Lemma \ref{lem: remove 1's} and Remark \ref{remark: remove rows}, we may remove the $\nu$th column and $i'$th row of the $t$-table, reducing $n, r$, and $m$ by one (by removing $\mu$ from $I$). If we prove that the theorem is true for $m + r -1 < n$, then we are able to prove it for all $m, r, n$. 
        \item If $m=2$, it will reduce to the first case. Set $t_0(i_0+1) = t_0(i_0 + 2) = t_1(i_1+1) = t_1(i_1+2)= 1$ (note that we required that $i_0 + 1 \neq i_1$ and $i_1 + 1 \neq i_0$, so these are 4 different values). 
        \item If $m > 2$, it will reduce to the second case. Since $m <  n-1$, there is some index $i' \in \{ 0, 1, \dots, n-1\}$ such that $i' \not \in I$. Let $j'$ be such that $i_{j'} < i' < i_{j'+1}$. This row will be where the `skip' in our step will be. \begin{enumerate}
            \item For $j \neq j'$, set $t_j(i_j+1) = t_j(i_{j+1}+1) = 1$ 
            \item Set $t_{j'}(i_{j'}+1) = t_{j'}(i_{j'}+2) = 1$ 
        \end{enumerate}  
        \item Notice that in steps 2 or 3, we would have placed a one in $4$ or $m+1$ columns respectively, and for each row, we have already placed 2 one's. We assumed that $m+r-1 < n$, so rearranging yields $r - 2 \leq  n - (m+2)$. Thus, we may pick $r-2$ columns that have not yet been filled and fill the column entirely with 1. 
        \item Set the rest of the $t_j(i)$ to $0$ 
    \end{enumerate}
    One may check that $(t_i)_i$ satisfies Condition \ref{condition}, and by removing columns of all 0's and 1's, the resulting chart will define one of the base cases above. 
\end{proof}

To illustrate the algorithm, consider the following examples, which reduce to the first and the second base case respectively. 

\begin{example}
    For the local model $M(3, 6, \{ 0, 4\})$, following the algorithm yields the $t$-table 
    \[
    \begin{tabular}{c|c|c|c|c|c|c}
        $i$ & 1 & 2 & 3 & 4 & 5 & 6 \\
        \hline
        $t_0$ & 1 & 1 & 1 & 0 & 0 & 0 \\
        $t_1$ & 0 & 0 & 1 & 0 & 1 & 1 
    \end{tabular}.
    \]
    We see that removing the 3rd and 4th column would yield the chart of the first base case. 
\end{example}

\begin{example}
    We notice that in order to apply step 1, we must have $n \geq 6$ and the reduction will result in the second base case. Let us reduce $M(4, 8, \{ 0, 1, 2, 4,6, 7\})$ to the second base case. For the first step, we may take $\mu = 7$ (and $\nu = 8$). For the third step, we may take $i' = 5$ and $j' = 3$. Thus, we obtain the chart 
    \[
    \begin{tabular}{c|c|c|c|c|c|c|c|c}
        $i$ & 1 & 2 & 3 & 4 & 5 & 6 & 7 & 8 \\
        \hline 
        $t_0$ & 1 & 1 & 0 & 1 & 0 & 0 & 0 & 1 \\
        $t_1$ & 0 & 1 & 1 & 1 & 0 & 0 & 0 & 1 \\
        $t_2$ & 0 & 0 & 1 & 1 & 1 & 0 & 0 & 1 \\
        $t_3$ & 0 & 0 & 0 & 1 & 1 & 1 & 0 & 1 \\
        $t_4$ & 1 & 0 & 0 & 1 & 0 & 0 & 1 & 1 \\
        $t_5$ & 1 & 1 & 0 & 1 & 0 & 0 & 0 & 1 \\
    \end{tabular}.
    \]
\end{example}

\printbibliography 

@online{hsy26,
      title={Cohen-Macaulayness of Local Models via Shellability of the Admissible Set}, 
      author={X. He and F. Schremmer and Q. Yu},
      date={2026},
      eprint={2603.05875},
      archivePrefix={arXiv},
      primaryClass={math.NT},
      url={https://arxiv.org/abs/2603.05875}, 
}

@book{matsumura,
    author = {H. Matsumura},
    translator = {Miles Reid},
    title = {Commutative Ring Theory},
    publisher = {Cambridge Univ. Press, Cambridge},
    collection={Cambridge Stud. Adv. Math.},
    year = {1987} 
}

@article{goertz,
    author = {U. Görtz},
    title = {On the flatness of models of certain Shimura varieties of PEL-type},
    journal = {Math. Ann.},
    volume = {321}, 
    year = {2001},
    pages = {689–-727}
}

@book{eisenbud,
    author = {D. Eisenbud},
    title = {Commutative Algebra},
    subtitle = {with a View Toward Algebraic Geometry},
    publisher = {Springer--Verlag, New York},
    year = {1995}, 
    collection = {Grad. Texts in Math.}
}

@book{RZ,
 author = {M. Rapoport and T. Zink},
 publisher = {Princeton Univ. Press, Princeton, NJ},
 title = {Period Spaces for $p$-divisible Groups},
 series = {Ann. of Math. Stud.}, 
 volume = {141}, 
 year = {1996}
}

@article{PRS-survey,
      title={Local models of Shimura varieties, I. Geometry and combinatorics}, 
      author={G. Pappas and M. Rapoport, and B. Smithling},
      year={2013},
      journal = {Handbook of Moduli},
      volume = {III.},
      pages = {135--219}
}

@article{Fakhruddin-Pilloni23, 
    title={Hecke Operators and the Coherent Cohomology of Shimura Varieties}, 
    volume={22},  
    number={1}, 
    journal={J. Inst. Math. Jussieu}, 
    author={N. Fakhruddin and V. Pilloni}, 
    year={2023}, 
    pages={1–-69}
}

@article{pappas-zhu,
    author = {G. Pappas and X. Zhu},
    title = {Local models of Shimura varieties and a conjecture of Kottwitz},
    journal = {Invent. Math.},
    year = {2013}, 
    volume={194}, 
    number={1}, 
    pages={147--254} 
}

@article{he13,
   title={Normality and Cohen–Macaulayness of local models of Shimura varieties},
   volume={162},
   number={13},
   journal={Duke Math. J.},
   author={X. He},
   year={2013},
   pages={2509 -- 2523}
}

@article{zhu14,
    title={On the coherence conjecture of Pappas and Rapoport},
    author={X. Zhu}, 
    volume={180}, 
    year={2014}, 
    pages={1--85}, 
    issue={1}, 
    journal = {Ann. of Math. (2)}
}

@article{hr22,
    author = {T. Haines and T. Richarz},
    title = {Normality and Cohen-Macaulayness of parahoric local models},
    journal = {J. Eur. Math. Soc.},
    year = {2022}, 
    volume={25}, 
    number={2}, 
    pages = {703--729}
}

@article{fh25,
    author = {N. Fakhruddin and T. Haines and J. Lourenco and T. Richarz},
    title = {Singularities of local models},
    journal = {Math. Ann.},
    year = {2025}, 
    pages={6205--6250}, 
    volume={391}
}

@article{bp00,
    author = {M. Brion and P. Polo},
    title = {Large Schubert varieties},
    journal = {Represent. Theory},
    year = {2000}, 
    pages = {97--126}, 
    number = {4}
}

@incollection{r90,
    author = {M. Rapoport},
    title = {On the Bad Reduction of Shimura Varieties},
    booktitle = {Automorphic Forms, Shimura Varieties, and L-functions, Vol. II},
    publisher = {Academic Press}, 
    year = {1990}, 
    editor = {L. Clozel and J. Milne}, 
    volume = {2}
}

@book{vakil,
    author = {Ravi Vakil},
    title = {The Rising Sea},
    subtitle = {Foundations of Algebraic Geometry},
    publisher = {Princeton Univ. Press},
    year = {2025}
}

@article{hpr20,
    author = {Xuhua He and Georgios Pappas and Michael Rapaport},
    title = {Good and semi-stable reduction of Shimura varieties},
    journal = {J. \'Ec. polytech. Math.},
    year = {2020}, 
    volume = {7}, 
    pages = {497--571}
}

@article{luo25,
    author = {Yu Luo},
    title = {On the moduli description of ramified unitary local models of signature $(n-1, 1)$},
    journal = {Math. Ann.}, 
    year = {2025}, 
    volume = {392},
    pages = {4661--4738}
}

@article{hls26, 
    author = {Qiao He and Yu Luo and Yousheng Shi},
    title = {The basic locus of ramified unitary Shimura varieties of signature $(n-1, 1)$ at maximal vertex level},
    date={2026},
      eprint={2502.06218},
      archivePrefix={arXiv},
      primaryClass={math.AG},
      url={https://arxiv.org/abs/2502.06218}, 
}

\end{document}